\documentclass[preprint]{imsart}

\usepackage{qzart}
\usepackage{amssymb}
\usepackage{multirow}
\usepackage{tabularx}

\begin{document}

\begin{frontmatter}

  \title{Two-Sample Instrumental Variable Analyses using Heterogeneous
    Samples}
  \runtitle{Two-Sample IV}

  \begin{aug}
    \author{Qingyuan Zhao\ead[label=e1]{qyzhao@wharton.upenn.edu}},
    \author{Jingshu Wang\ead[label=e2]{jingshuw@wharton.upenn.edu}},
    \author{Wes Spiller\ead[label=e3]{wes.spiller@bristol.ac.uk}}
    \author{Jack Bowden\ead[label=e4]{jack.bowden@bristol.ac.uk}},
    \and
    \author{Dylan S.\ Small\ead[label=e5]{dsmall@wharton.upenn.edu}}

    \runauthor{Q. Zhao et al.}

    \affiliation{University of Pennsylvania and University of Bristol}

    \address{Department of Statistics, The Wharton School, University
      of Pennsylvania, USA \printead{e1,e2,e5}.}

    \address{MRC Integrative Epidemiology Unit, University of Bristol,
      UK \printead{e3,e4}.}

  \end{aug}

  \begin{abstract}
    Instrumental variable analysis is a widely used method to estimate causal
    effects in the presence of unmeasured
    confounding. When the instruments, exposure and outcome are not
    measured in the same sample, \citet{angrist1992effect} suggested to
    use two-sample instrumental variable (TSIV) estimators that use sample
    moments from an instrument-exposure sample and an instrument-outcome
    sample. However, this method is biased if the two samples
    are from heterogeneous populations so that the distributions of the
    instruments are different. In linear structural equation models, we
    derive a new class of TSIV estimators
    that are robust to heterogeneous samples under the key assumption that
    the structural relations in the two samples are the same. The widely
    used two-sample two-stage least squares estimator belongs to this
    class. It is generally not asymptotically efficient, although we
    find that it performs similarly to the optimal TSIV estimator in
    most practical situations. We then attempt to relax the linearity
    assumption. We find that, unlike one-sample analyses, the TSIV
    estimator is not robust to misspecified exposure
    model. Additionally, to nonparametrically identify the magnitude
    of the causal effect, the noise in the exposure must
    have the same distributions in the two samples. However, this
    assumption is in general untestable because the exposure is not
    observed in one sample. Nonetheless, we may still identify the sign of
    the causal effect in the absence of homogeneity of the noise.
  \end{abstract}


\end{frontmatter}

\section{Introduction}
\label{sec:introduction}

When randomized controlled experiments are not feasible, instrumental
variable (IV) analysis is a widely used method to estimate causal
effect in the presence of unmeasured confounding. A typical instrumental
variable estimator such as the two-stage least squares (TSLS) uses
sample moments (e.g.\ covariance matrices) of the instrument-exposure
relationship and the instrument-outcome relationship. In an
influential article, \citet{angrist1992effect} noticed that the two
sets of moments can indeed be estimated from different samples, though
this idea can be dated back to at least
\citet{klevmarken1982missing}. This
method, often referred to as the two-sample instrumental variable
(TSIV) estimator, is frequently used in econometrics
\citep{inoue2010two}.

One of the most exciting recent applications of IV analysis is
in genetic epidemiology where genetic variants are used as the
instruments
\citep{davey2003mendelian,lawlor2008mendelian,burgess2015review}. This
method is known as ``Mendelian randomization'' to epidemiologists,
because the genotypes are governed by Mendel's Second Law of
independent assortment and thus have a strong rationale for being
independent of common postnatal source of
confounding. More recently, there has been growing interest in using
two-sample Mendelian randomization that take advantage of large existing
Genome-Wide Association Studies (GWAS), as it is often easier to
find two GWAS in which one measures the genotypes
and the exposure and the other one measures the genotypes and the
disease than to find a single GWAS that measures all three types of
variables
\citep{pierce2013efficient,davey2014mendelian,burgess2015using,gamazon2015gene,lawlor2016commentary}.

Since Mendelian randomization is a special case of instrumental
variable analysis in which genetic variants are used as instruments,
one would expect that two-sample Mendelian randomization is merely
a different application of the existing TSIV estimators. However,
there is a subtle but important difference between two-sample Mendelian
randomization and the existing applications of TSIV in economic applications. To the
best of our knowledge, with the exception of
\citet{graham2016efficient} who considered a general data combination
problem including just-identified TSIV, all the TSIV estimators
previously proposed in econometrics assumed that the two datasets are
sampled from the same population
\citep{angrist1992effect,ridder2007econometrics,inoue2010two,pacini2016robust}. This
is
usually not a problem in the economic applications using time-invariant
instrumental variables \citep{jappelli1998testing} such as quarter of birth
\citep{angrist1992effect} and sex composition of the children in the
household \citep{currie2000public}. However, this
assumption does not hold in two-sample Mendelian
randomization, as the two GWAS usually consist of different cohort
studies and thus represent different populations. \Cref{tab:tsmr}
shows an example of two-sample Mendelian randomization in which the
distribution of the genetic instruments are clearly different in the
different populations.


\renewcommand{\arraystretch}{1.25}

\begin{table}[t]
  \centering
  \caption{Heterogeneous distribution of genetic instruments in
    different populations. The minor allele frequencies in
    \Cref{tab:tsmr-maf} and link disequilibrium $r^2$ in
    \Cref{tab:tsmr-ld} are obtained from the 1000 Genome Project
    available in online databases dbSNP \citep{sherry2001dbsnp} and
    LDlink \citep{machiela2015ldlink}. The SNPs are selected from the
    real data analysis in \Cref{sec:real-data-analysis}.}
  \label{tab:tsmr}
  \begin{subtable}{\textwidth} \centering
    \caption{Minor allele frequencies in different populations.}
    \label{tab:tsmr-maf}
    \begin{tabular}{ll|ccccc}
      \hline
      &  & \multicolumn{5}{c}{Minor allele frequency} \\
      &  & African & East Asian & European & South
      Asian
      &
      American
      \\
      & & ($n = 1322$) & ($n = 1008$) & ($n = 1006$) & ($n = 978$) &
      ($n
      =
      694$) \\
      \hline
      &rs13021737  & 0.095 & 0.087 & 0.174 & 0.13 & 0.14 \\
      SNP &rs1421085  & 0.056 & 0.169 & 0.432 & 0.31 & 0.24 \\
      &rs6567160  & 0.220 & 0.183 & 0.240 & 0.32 & 0.13 \\
      \hline
    \end{tabular}
  \end{subtable}

  \begin{subtable}{\textwidth} \centering
    \caption{Linkage disequilibrium (measured by $r^2$, the square of
      the correlation coefficient of allele indicators) in different
      populations.}
    \label{tab:tsmr-ld}
    \begin{tabular}{ll|ccccc}
      \hline
      &  & \multicolumn{5}{c}{Linkage disequilibrium ($r^2$)} \\
      &  & African & East Asian & European & South
      Asian
      &
      American
      \\
      & & ($n = 1322$) & ($n = 1008$) & ($n = 1006$) & ($n = 978$) &
      ($n
      =
      694$) \\
      \hline
      &(rs13021737, rs6731348)  & 0.378 & 1.0 & 0.993 & 0.865 & 0.965 \\
      SNP pair & (rs13021737, rs4854344) & 0.917 & 1.0 & 0.993 & 0.865 & 0.988 \\
      &(rs6731348, rs4854344)  & 0.387 & 1.0 & 0.986 & 1.0 & 0.953 \\
      \hline
    \end{tabular}
  \end{subtable}

\end{table}

The goal of this paper is to clarify the consequences of heterogeneous
samples to the identification, estimation, and robustness of TSIV
analyses. After setting up the TSIV problem and reviewing the
literature (\Cref{sec:backgr-tsiv-analys}), we will derive a new class of
TSIV estimators using the generalized method of moments (GMM) that can
utilize two heterogeneous samples under a linear IV model
(\Cref{sec:tsiv-estim-using}). The commonly-used two-sample two-stage
least squares (TSTSLS) belongs to this class of estimators, but unlike
the case with homogeneous samples, it is no longer the most efficient
estimator in this class. Another
interesting question raised by epidemiologists and geneticists is how
far we can get by using just public summary statistics of GWAS
\citep{lawlor2016commentary,barbeira2016metaxcan}. Our calculations
show that, to use correlated genetic IVs without individual-level
data, it is necessary to use their covariance matrices (in both
samples) to compute any TSIV estimator and its asymptotic
variance. Unfortunately,
the covariance information is often unavailable in the current GWAS
summary databases, though it is possible to approximate the covariance
matrices using external datasets such as the \citet{10002015global}.

We will then turn to relax the linearity assumption in
\Cref{sec:relax-invar-assumpt,sec:relax-line-instr,sec:relax-line-expos}. Compared
to the same problem in the one-sample or the homogeneous two-sample
setting, a key distinction is
that we also need the structural relationships between the IV and the
exposure and the distributions of the noise variables to be invariant
in the two samples. Unfortunately, these assumptions are untestable
using empirical data because we do not observe the exposure in both
samples. In the absence of these assumptions, we show that one may
still identify the sign of the causal effect.

Next we will use simulations to study the numerical properties of the TSIV estimators (\Cref{sec:simulation}). We find that although the
asymptotic efficiency of TSTSLS is suboptimal theoretically, the
difference in practice is most of the time minuscule. We will also
examine the bias of the TSIV estimators when the
instrument-exposure equation is misspecified or the ``homogeneous
noise'' assumption is violated. We will also compare the results of
the TSIV analyses with the classical one-sample analyses using a
real Mendelian randomization dataset (\Cref{sec:real-data-analysis}).
Finally, we will summarize the theoretical and empirical findings in
\Cref{sec:summary-discussion}. Although we
will be using Mendelian randomization as the motivating application in
the investigation below,
we expect the statistical methods, identification results and high-level
conclusions in this paper can be applied to TSIV
analyses in other fields as well.

\section{Background on TSIV analyses}
\label{sec:backgr-tsiv-analys}

In this section we set up the TSIV problem and review the related
literature. For simplicity of exposition, throughout the paper we
consider only one endogenous exposure variable and no other
exogenous covariates for adjustment.
Most of our derivations can be easily generalized to the case
of multiple endogenous variables and multiple exogenous covariates.

\subsection{Problem setup}
\label{sec:setup}

We begin by introducing some notational conventions. We use lower-case
letters, bold lower-case letters, bold
upper-case letters, and Greek letters to indicate, respectively,
deterministic or random scalars, vectors, matrices, and parameters in
the model. Superscripts $s$, $a$, $b$ are reserved to indicate the
sample. Subscripts are used to index the observations in each sample.

Suppose we have independent samples
$(\mathbf{z}_i^s,x_i^s,y_i^s),~i=1,2,\dotsc,n^s$,
from two populations, $s = a$ and $s = b$, where $\mathbf{z} \in \mathbb{R}^q$ is
a vector of instrumental variables, $x$ is the exposure variable,
and $y$ is the outcome variable. More compactly, we can write the data
in each sample as a matrix $\mathbf{Z}^s \in \mathbb{R}^{n^s \times
  q}$ and two vectors $\mathbf{x}^s,\mathbf{y}^s \in
\mathbb{R}^{n^s}$. Next we describe the general setting in this
paper.

\begin{assumption} \label{assump:npsem}
  The data are generated from the following
  nonparametric structural equation model (SEM). For $s \in \{a,b\}$,
  \begin{align}
    y_i^s &= g^s(x_i^s,\mathbf{z}_i^s,u_i^s), \label{eq:npsem-outcome} \\
    x_i^s &= f^s(\mathbf{z}_i^s, v_i^s), \label{eq:npsem-exposure}
  \end{align}
  where the functions $g^s$, $f^s$ are unknown and the random variables
  $(u^s_i,v^s_i,\mathbf{z}^s_i),~i=1,\dotsc,n^s$ are independent and
  identically distributed within each sample.
\end{assumption}

Hereafter, \eqref{eq:npsem-outcome} will be called the \emph{exposure-outcome
  equation} or simply the \emph{outcome equation}, and \eqref{eq:npsem-exposure} the \emph{instrument-exposure
  equation} or the \emph{exposure} equation. The exposure variable $x$
is called \emph{endogenous} if
$v \not \independent u$ (so $x \not \independent u$). In this case, a
plain regression of $y$ on $x$ would lead to biased estimate of the
causal effect of $x$.

There are three necessary conditions for $\mathbf{z}$ to be
\emph{valid} instrumental variables:
$\mathbf{z}$ must be correlated with $x$, $\mathbf{z}$ must be
independent of the unmeasured confounder(s), and $\mathbf{z}$ must affect
the outcome $y$ only through $x$ (exclusion restriction). These
assumptions are usually stated in the potential outcome language
\citep{angrist1996identification}. Translating these into structural
equation models, we need to assume the following core IV assumptions:
\begin{assumption} \label{assump:valid-iv}
  (Validity of IV) For $s \in \{a,b\}$, the exposure equation $f^s$ is
  a non-constant function on the support of $\mathbf{z}^s$, $\mathbf{z}_i^s \independent
  (u_i^s,v_i^s)$, and the outcome equation $g^s$ does not depend on
  $\mathbf{z}^s$.
\end{assumption}

Next we describe the one-sample and two-sample IV problems:

{\noindent \bf The classical IV problem:} Suppose we observe
$\mathbf{Z}^a$, $\mathbf{x}^a$, and $\mathbf{y}^a$ in the first
sample. If $x$ is endogenous, what can we learn about the outcome
equation \eqref{eq:npsem-outcome} (how $g^a$ behaves as a function of
$x^a$, a.k.a. the ``causal effect'' of $x$ on $y$) by using the
instrumental variables $\mathbf{Z}^{a}$?

{\noindent \bf The two-sample IV problem:} Suppose
only $\mathbf{Z}^a$, $\mathbf{x}^a$, $\mathbf{Z}^b$, and
$\mathbf{y}^b$ are observed (in other words $\mathbf{y}^a$ and
$\mathbf{x}^b$ are not observed). If $x$ is endogenous, what can we
learn about the outcome equation \eqref{eq:npsem-outcome}?

In the classical one-sample setting, the valid IV assumption
(\Cref{assump:valid-iv}) is not sufficient to identify the
causal effect of $x$ on $y$. Further assumptions are required to
identify the causal effect. The simplest and most widely studied
setting is when the instrument-exposure and exposure-outcome equations
are both linear (linearity of the exposure equation is not necessary,
see \Cref{sec:relax-line-instr}):
\begin{assumption} \label{assump:linearity}
  (Linearity) For $s \in \{a,b\}$, (\ref{assump:linearity}-1) $g^s(x_i,u_i) = \beta^s x_i +
  u_i$; and (\ref{assump:linearity}-2) $f^s(\mathbf{z}_i,v_i) = (\bm{\gamma}^s)^T \mathbf{z}_i + v_i$.
\end{assumption}
Under \Cref{assump:linearity}, the structural equations
\eqref{eq:npsem-outcome} and \eqref{eq:npsem-exposure} can be written in a more
compact form: for $s\in\{a,b\}$,
\begin{equation}
  \label{eq:model-simple}
  \begin{split}
    &\mathbf{y}^s = \mathbf{x}^s \beta^s + \mathbf{u}^s, \\
    &\mathbf{x}^s = \mathbf{Z}^s \bm{\gamma}^s + \mathbf{v}^s.
  \end{split}
\end{equation}
Without loss of generality, we assume the expected values of
$\mathbf{z}$, $u$ and $v$ in both samples are $0$. Otherwise we can just
add intercept terms to \eqref{eq:model-simple}.

Another commonly used assumption is monotonicity which leads to the
identification of the local average treatment effect (LATE), see
\Cref{assump:monotonicity} in \Cref{sec:relax-line-expos}. We will see
that in the two-sample setting, even more assumptions are needed to
identify the causal effect.

\subsection{Literature review}
\label{sec:literature-review}

Next we give a literature review on instrumental variables
regression. Our goal is to not give the most comprehensive review of
this massive literature, but rather to outline some key ideas to aid
us in the investigation of the TSIV estimators using heterogeneous
samples. We will also discuss problems (such as weak IV bias
and invalid IV bias) that are commonly encountered
in Mendelian randomization studies.

\subsubsection{One-sample IV estimators.}
\label{sec:class-iv-estim}

IV methods were developed in early research on structural/simultaneous
equation modeling by
\citet{wright1928tariff,anderson1949estimation,theil1958economic}
among many others. For simplicity, when
considering the one-sample IV problem below we shall ignore the
superscript $a$.  The most important and widely used estimator in the
classical setting is the two-stage least squares (TSLS), where the
exposure $x$ is first regressed on the IVs $\bm z$ (first-stage
regression) using least squares and the outcome $y$ is then regressed
on the predicted exposure from the first-stage regression using
another least squares. The TSLS estimator can be concisely written
using the projection matrix $\mathbf{P}_{\mathbf z} = \mathbf{Z} (\mathbf{Z}^T
\mathbf{Z})^{-1} \mathbf{Z}^T$:
\[
\hat{\beta}_{\text{TSLS}} = (\mathbf{x}^T \mathbf{P}_{\mathbf z}
\mathbf{x})^{-1} (\mathbf{x}^T \mathbf{P}_{\mathbf z} \mathbf{y}).
\]
Other classical IV estimators include the limited information maximum
likelihood (LIML) \citep{anderson1949estimation} and \citet{fuller1977some}'s modified LIML estimator. All these estimators
belong to the general $K$-class estimators \citep{theil1958economic}. For
a more comprehensive textbook treatment of the classical IV estimators, we
refer the reader to \citet{davidson1993estimation}.

There is also considerable effort to relax the homogeneous causal
effect assumption in \eqref{eq:model-simple}. The most influential
approach is the LATE framework \citep{imbens1994identification,baker1994paired,angrist1996identification} that
will be discussed in detail in \Cref{sec:relax-line-expos}. See
\citet{abadie2003semiparametric,ogburn2015doubly} for some recent
methodological developments in this direction. Another
approach is to assume all the effect modifiers in the exposure- and
outcome-equations are observed
\citep{hernan2006instruments,wang2018bounded}. \citet{baiocchi2014instrumental}
gives a comprehensive review of one-sample IV estimators in biomedical applications.

\subsubsection{Two-sample IV estimators.}
\label{sec:two-sample-instr}

The idea of using different samples to estimate moments can be dated
back to \citet{klevmarken1982missing} and this proposal becomes
popular in econometrics after \citet{angrist1992effect}. In a later
article, \citet{angrist1995split} further argued to routinely use the
split-sample TSLS estimator so that weak instrument biases the estimator
towards $0$ instead of towards the ordinary least squares (OLS)
estimator. \citet{inoue2010two} compared the asymptotic distributions
of alternative TSIV estimators. They found that the TSTSLS estimator
is not only  more efficient than the covariance-based TSIV estimator,
but also achieves asymptotic efficiency in the class of limited
information estimators. \citet{ridder2007econometrics} considered a
more general form of TSIV estimator and derived its asymptotic
distribution. More recently, \citet{pacini2016robust} derived
heteroskedasticity-robust variance estimator of TSTSLS and
\citet{pacini2018two} derived a semiparametrically efficient TSIV
estimator with interval-censored covariates. All the references above
considered the TSIV problem with homogeneous samples. The only
exceptions we know are \citet{graham2016efficient} who considered a
general data combination problem which includes the just-identified
TSIV, and a working version of \citet{inoue2010two} who considered
different sampling rates dependent on the instruments.

\subsubsection{Summary-data Mendelian randomization.}
\label{sec:summ-data-mend}

Since Mendelian randomization is just a special case of IV analyses
where genetic variation is used as the IV, all the one-sample or
two-sample methods mentioned above can be directly
applied. However, when conducting Mendelian randomization studies we
only have access to ``summary data'' that only contain the marginal
regression coefficients and their standard errors. For example, let
the estimated regression coefficient of $\bm y$ on $\bm Z_{\cdot j}$ be
$\hat{\Gamma}_j$ and the coefficient of $\bm x$ on $\bm Z_{\cdot j}$ be
$\hat{\gamma}_j$. Then \citet{wald1940fitting}'s ratio estimator of
the causal effect using the $j$-th instrument is given by
$\hat{\beta}_j = \hat{\Gamma}_j / \hat{\gamma}_j$. This is equivalent
to using a single instrument in TSLS. The statistical problem is then
to combine the individual estimators, like in a meta-analysis, to
produce a single efficient and robust estimator.

The above summary-data Mendelian randomization design has wide
applicability in practice \citep{burgess2015using} and there is a lot
of ongoing efforts in developing public databases and software
platforms \citep{hemani2018mr}. In human genetics, Mendelian
randomization is used as a tool for gene testing and discovery
\citep{gamazon2015gene}. On the methodological side, the
commonly used meta-analysis estimators in this problem include Egger
regression \citep{bowden2015mendelian} and weighted median
\citep{bowden2016consistent}. More recently,
\citet{zhao2018statistical} proposed to treat summary-data Mendelian
randomization as a errors-in-variables regression problem to develop
more efficient and robust estimators.

\subsubsection{Weak IVs and invalid IVs.}
\label{sec:weak-instruments}

Finally we want to briefly mention a critical problem that plagues
many IV analyses---invalidity of the
instruments. One such problem is the weak instrument bias that occurs
when the IVs $\bm z$ are only weakly associated with the exposure
$x$. In this case, the classical IV estimators are usually biased
towards the OLS estimator in one-sample setting or towards 0 in the
two sample setting. This problem has been well studied in the
one-sample setting, see \citet{stock2002survey} for a comprehensive
survey. In Mendelian randomization it is common to have many weak
instruments. In this regime, LIML-like estimators are asymptotically
unbiased but the asymptotic variance needs to be carefully derived
\citep{hansen2008estimation}. More recently, \citet{choi2018weak}
studied this problem in the two-sample setting and
\citet{zhao2018statistical} proposed robust statistical inference in
summary-data Mendelian randomization with many weak instruments.

Compared to weak IV bias, more serious problems can be caused by invalid
instruments that are dependent on unmeasured confounders or violate
the exclusion restriction assumption. In classical IV analyses with
one or just a few IVs, the analyst must use domain knowledge to
justify the validity of the instruments. In Mendelian randomization,
the exclusion restriction assumption may be violated due to a
genetic phenomenon called pleiotropy
\citep{davey2003mendelian}. Fortunately we often have dozens and
hundreds of independent genetic instruments, and it is possible to use additional
assumptions such as sparsity of invalid IVs
\citep{kang2016instrumental} or balanced direct effects
\citep{bowden2015mendelian} to identify and estimate the causal
effect.

For the rest of this paper, we will assume all the IVs are strong and
valid. Our goal is to show that, in addition to the weak and invalid
IV problems mentioned above, heterogeneity of the samples can bring new
challenges to the inference and interpretation of TSIV analyses.

\section{Linear TSIV estimators using heterogeneous samples}
\label{sec:tsiv-estim-using}

In \Cref{assump:npsem,assump:valid-iv,assump:linearity}, we have been
stating our assumptions separately for each sample. If the
structural relationships can be arbitrarily different in the two
samples, it is obviously hopeless to solve the endogeneity problem
with two partially observed samples. We use the next two assumptions to
link the structural equations in the two samples.

\begin{assumption} \label{assump:invariance}
  (Structural invariance) $g^a = g^b$, $f^a = f^b$.
\end{assumption}

\begin{assumption} \label{assump:invariance-noise}
  (Sampling homogeneity of the noise variables) $(u_i^a,v_i^a) \overset{d}{=}
  (u_j^b,v_j^b)$ for any $i=1,\dotsc,n^a$, $j=1,\dotsc,n^b$.
\end{assumption}

Both assumptions put restrictions on the heterogeneity of the two
samples. To distinguish structural and distributional assumptions,
we use different words---``invariance'' and ``homogeneity''---to refer
to these assumptions. Under
\Cref{assump:invariance,assump:invariance-noise}, the only
heterogeneity between the two samples comes from the distribution of
the instruments.
In linear SEMs, \Cref{assump:invariance-noise} is not
required (see \Cref{sec:relax-invar-assumpt}), but it is generally
necessary in nonparametric SEMS because we do not specify the forms of the
functions $f$ and $g$ in \Cref{assump:npsem,assump:invariance}.

In this Section we will study TSIV estimators in the linear SEM
\eqref{eq:model-simple}. In this case, structural invariance
or \Cref{assump:invariance} implies
that $\beta^a = \beta^b = \beta$, $\bm{\gamma}^a = \bm{\gamma}^b =
\bm{\gamma}$. Our inferential target is the parameter
$\beta$, which is interpreted as the causal effect of $x$ on $y$.

We introduce some notations for the covariance parameters in this
model. For $s \in \{a,b\}$, denote the population covariances as
$\mathrm{Cov}(\mathbf{Z}^s) = \bm{\Sigma}^s_{zz}$,
$\mathrm{Cov}(\mathbf{Z}^s,\mathbf{x}^s) = \bm{\Sigma}_{zx}^s$,
$\mathrm{Var}(\mathbf{u}^s) = (\sigma_{uu}^s)^2$,
$\mathrm{Var}(\mathbf{v}^s) = (\sigma_{vv}^s)^2$,
$\mathrm{Cov}(\mathbf{u}^s,\mathbf{v}^s) =
\sigma_{uv}^s$. Denote the sample covariance matrices as (recall that
we assume all the random variables have mean $0$)
\[
\mathbf{S}^s_{zy} = (\mathbf{Z}^s)^T\mathbf{y}^s/n^s,~\mathbf{S}^s_{zx} =
(\mathbf{Z}^s)^T \mathbf{x}^s/n^s,~\mathbf{S}^s_{zz} = (\mathbf{Z}^s)^T
\mathbf{Z}^s/n^s.
\]

We use the generalized method of moments (GMM) in
\citet{hansen1982large} to estimate $\beta$ under
\Cref{assump:npsem,assump:valid-iv,assump:invariance,assump:invariance-noise,assump:linearity}. Consider
the following moment function of $\beta$:
\[
\mathbf{m}_n(\beta) = (\mathbf{S}^b_{zz})^{-1} \mathbf{S}^b_{zy} -
(\mathbf{S}^a_{zz})^{-1} \mathbf{S}^a_{zx} \beta.
\]
Compared to the moment function defined in \citet{angrist1992effect},
we added the normalization terms $(\mathbf{S}^a_{zz})^{-1}$ and
$(\mathbf{S}^b_{zz})^{-1}$ because $\bm{\Sigma}_{zz}^a$ and
$\bm{\Sigma}_{zz}^b$ can be different in the heterogeneous two-sample
setting. To differentiate between an arbitrary value of $\beta$ and
the true value of $\beta$, we use $\beta_0$ to denote the true value
in this section. First, we check the moment conditions $\mathbb{E}[\mathbf{m}_n(\beta)] =
\mathbf{0}$ identifies $\beta_0$ by showing
$\mathbf{m}_n(\beta_0) \overset{d}{\to} \mathbf{0}$. To see this, notice
that
\begin{equation} \label{eq:gn-mean}
  \begin{split}
    \mathbf{m}_n(\beta) &= (\mathbf{S}^b_{zz})^{-1} (\mathbf{Z}^b)^T(\mathbf{Z}^b
    \bm{\gamma} \beta_0 + \mathbf{v}^b \beta_0 + \mathbf{u}^b) / n_b -
    (\mathbf{S}^a_{zz})^{-1} (\mathbf{Z}^a)^T (\mathbf{Z}^a \bm{\gamma} +
    \mathbf{v}^a) \beta /n_a \\
    &= \bm \gamma (\beta_0 - \beta) + (\mathbf{S}^b_{zz})^{-1} (\mathbf{Z}^b)^T(\mathbf{v}^b \beta +
    \mathbf{u}^b) / n_b - (\mathbf{S}^a_{zz})^{-1} (\mathbf{Z}^a)^T
    \mathbf{v}^a \beta /n_a.
  \end{split}
\end{equation}
It is easy to see that $\mathbf{m}_n(\beta_0)$ has mean $0$ and converges to $0$
in probability. The key in \eqref{eq:gn-mean} is that the normalization by
$(\mathbf{S}^a_{zz})^{-1}$ and $(\mathbf{S}^b_{zz})^{-1}$ makes sure
the first term on the right hand side is $0$ when $\beta = \beta_0$.

Next, let $\mathbf{W} \in \mathbb{R}^{q \times q}$ be a positive definite
weighting matrix. The class of TSIV
estimators of $\beta$ is given by
\begin{equation} \label{eq:beta-hat}
  \begin{split}
    \hat{\beta}_{n,\mathbf{W}} &= \argmin_{\beta}\,\mathbf{m}_n(\beta)^T \mathbf{W}
    \mathbf{m}_n(\beta) \\
    &= [(\mathbf{S}^a_{zx})^T (\mathbf{S}^a_{zz})^{-1}
    \mathbf{W} (\mathbf{S}^a_{zz})^{-1} \mathbf{S}^a_{zx}]^{-1}
    [(\mathbf{S}^a_{zx})^T (\mathbf{S}^a_{zz})^{-1}
    \mathbf{W} (\mathbf{S}^b_{zz})^{-1} \mathbf{S}^b_{zy}].
  \end{split}
\end{equation}
Using the general theory for GMM \citep{hansen1982large}\footnote{As
  pointed out by a reviewer, our application of the GMM theory is a
  bit non-standard because GMM usually starts with moment functions
  that only depend on \emph{one} data point. Nevertheless, the asymptotic
  normality still goes through by a similar ``sandwich'' argument
  because $m_n(\beta)$ is asymptotically normal.}, the
asymptotic variance of $\hat{\beta}_{n,\mathbf{W}}$ is given by
\begin{equation} \label{eq:beta-hat-variance}
  \begin{split}
    &\mathrm{Var}(\hat{\beta}_{n,\mathbf{W}}) \approx \\
    &[(\mathbf{S}^a_{zx})^T (\mathbf{S}^a_{zz})^{-1} \mathbf{W}
    (\mathbf{S}^a_{zz})^{-1} \mathbf{S}^a_{zx}]^{-1}
    (\mathbf{S}^a_{zx})^T (\mathbf{S}^a_{zz})^{-1} \mathbf{W} \bm{\Omega}_n
    \mathbf{W} (\mathbf{S}^a_{zz})^{-1} \mathbf{S}^a_{zx}
    [(\mathbf{S}^a_{zx})^T (\mathbf{S}^a_{zz})^{-1} \mathbf{W}
    (\mathbf{S}^a_{zz})^{-1} \mathbf{S}^a_{zx}]^{-1},
  \end{split}
\end{equation}
where $\bm{\Omega}_n$ is the variance of
$\mathbf{m}_n(\beta_0)$.

The optimal $\mathbf{W}$ in this class of estimators is given by $\mathbf{W} \propto
\bm{\Omega}_n^{-1}$. Next we compute $\bm{\Omega}_n$. It is easy to
see that
\begin{align}
  \mathrm{Var}(\mathbf{m}_n(\beta_0)\,|\,\mathbf{Z}^a,\mathbf{Z}^b)
  =& \mathrm{Var}((\mathbf{S}^b_{zz})^{-1} \mathbf{S}^b_{zy}|\mathbf{Z}^b) +
  \mathrm{Var}((\mathbf{S}^a_{zz})^{-1}
  \mathbf{S}^a_{zx}\beta_0|\mathbf{Z}^a) \nonumber \\
  =&\frac{1}{n_b} (\mathbf{S}_{zz}^b)^{-1} [\beta_0^2 (\sigma_v^b)^2 + 2 \beta_0
  \sigma_{uv}^b + (\sigma_u^b)^2] + \frac{1}{n_a}
  (\mathbf{S}_{zz}^a)^{-1} [\beta_0^2 (\sigma_v^a)^2] \nonumber \\
  =& \frac{1}{n_b} (\mathbf{S}_{zz}^b)^{-1} \mathrm{Var}(y_i^b|\mathbf{z}_i^b) +
  \frac{1}{n_a}(\mathbf{S}_{zz}^a)^{-1} \beta_0^2
  \mathrm{Var}(x_i^a|\mathbf{z}_i^a). \label{eq:omega-2}
\end{align}
In other words, the conditional variance of $\mathbf{m}_n(\beta_0)$ is
the sum of the variance of the coefficient of the outcome-instrument
regression (in sample $b$) and $\beta_0^2$ times the variance of the
coefficient of the exposure-instrument regression (in sample
$a$). 
Equation \eqref{eq:omega-2} means that to estimate $\mathbf{\Omega}_n$ and the variance of
$\hat{\beta}_{n,\mathbf{W}}$ for any given $\mathbf{W}$, we just need
to estimate the noise variances of the outcome-instrument and
exposure-instrument regressions. Weak instrument bias may occur when
the magnitude of $\bm \gamma$ is small comparing to $\sigma_v^2$. In
this case the asymptotics presented here may be inaccurate and the
TSIV estimators are biased towards 0.


The asymptotically efficient two-sample IV estimator is
$\hat{\beta}_{n,\hat{\bm{\Omega}}_n^{-1}}$. Its asymptotic
variance is given by
\begin{equation} \label{eq:tsiv-var}
  \mathrm{Var}(\hat{\beta}_{n,\hat{\bm{\Omega}}_n^{-1}}Z) \approx
  [(\bm{\Sigma}^a_{zx})^T (\bm{\Sigma}^a_{zz})^{-1} \bm{\Omega}^{-1}_n
  (\bm{\Sigma}^a_{zz})^{-1} \bm{\Sigma}^a_{zx}]^{-1},
\end{equation}
which can be consistently estimated by $[(\hat{\bm{\gamma}}^a)^{-1}
\hat{\bm{\Omega}}_n^{-1} \hat{\bm{\gamma}}^a]^{-1}$.

We would like to make five remarks on the new class of TSIV estimators.

\begin{remark}
  When the weighting matrix $\bm{W}$ is chosen as $\mathbf{S}_{zz}^b$, the estimator reduces to the
  two-sample two-stage least squares (TSTSLS) estimator. To see this, let $\hat{\bm{\gamma}} = (\mathbf{S}^a_{zz})^{-1}
  \mathbf{S}^a_{zx}$ and $\hat{\mathbf{x}}^b = \mathbf{Z}^b
  \hat{\bm{\gamma}}$ be the predicted values. Then the TSTSLS estimator is
  defined as
  \[
  \hat{\beta}_{\text{TSTSLS}} = [(\hat{\mathbf{x}}^b)^T \hat{\mathbf{x}}^b]^{-1}
  (\hat{\mathbf{x}}^b)^T \mathbf{y}^b.
  \]
  It is easy to verify that $\hat{\beta}_{n,\mathbf{S}_{zz}^b} =
  \hat{\beta}_{\text{TSTSLS}}$. Thus unlike in the classical one-sample and
  homogeneous two-sample settings, TSTSLS is generally not efficient
  in the class of linear TSIV estimators when the two-samples are
  heterogeneous. To the best of our knowledge, this results is not
  known previously. Also, notice that the conventional covariance
  estimator based on sample covariance matrices is generally
  biased. In the exact-identified case ($q = 1$), the two-sample
  covariance estimator is
  \[
  \hat{\beta}_{\text{TSCOV}}= (\mathbf{s}^a_{zx})^{-1} \mathbf{s}^b_{zy}
  \overset{p}{\to} \beta_0 \cdot (\sigma^b_{zz} / \sigma^a_{zz}).
  \]
  In the homogeneous TSIV problem, the TSCOV estimator is not biased
  but less efficient than TSTSLS \citep{inoue2010two}. The
  inconsistency of TSCOV in heterogeneous TSIV problem is also noticed
  in \citet[footnote 1]{inoue2010two}.
\end{remark}

\begin{remark} \label{rem:efficiency}
  Notice that $\bm{\Omega}_n$ is a weighted sum of
  $(\bm{\Sigma}_{zz}^a)^{-1}$ and
  $(\bm{\Sigma}_{zz}^b)^{-1}$. In the homogeneous TSIV problem where $\bm{\Sigma}_{zz}^a =
  \bm{\Sigma}_{zz}^b$, we have $\bm{\Sigma}_{zz}^b
  \propto \bm{\Omega}_n^{-1}$ and hence the TSTSLS estimator is
  efficient in the class of TSIV estimator \eqref{eq:beta-hat}. This
  is consistent with the conclusions of \citet[Theorem
  1]{inoue2010two}. In general, the efficiency of the TSTSLS
  estimator (relative to the most efficient TSIV estimator
  $\hat{\beta}_{n,\hat{\bm{\Omega}}_n^{-1}}$) depends on the
  difference between $(\bm{\Sigma}_{zz}^a)^{-1}$ and
  $(\bm{\Sigma}_{zz}^b)^{-1}$, the ratio of $n^a$ and $n^b$, and the
  ratio of $\mathrm{Var}(y_i^b|\mathbf{z}_i^b)$ and
  $\beta_0^2\mathrm{Var}(x_i^b|\mathbf{x}_i^b)$. In
  most cases we expect the covariance
  structures of the instrumental variables are not too different in
  the two samples and the last ratio to be not too small, so the
  TSTSLS estimator has great relative efficiency. We will see that
  TSTSLS and the optimal TSIV estimator have very similar performance
  in simulations (\Cref{sec:simulation}).
\end{remark}

\begin{remark}
  A naive estimator of the asymptotic variance of $\hat{\beta}_{\text{TSTSLS}}$
  is simply the variance of the coefficient in the second-stage regression:
  \[
  \hat{\sigma}^2_{\text{naive}}(\hat{\beta}_{\text{TSTSLS}}) = [(\hat{\mathbf{x}}^b)^T
  \hat{\mathbf{x}}^b]^{-1} \widehat{\mathrm{Var}}(y_i^b|\hat{x}^b_i) \to
  [(\bm{\Sigma}^a_{zx})^T (\bm{\Sigma}^a_{zz})^{-1} \tilde{\bm{\Omega}}_n^{-1}
  (\bm{\Sigma}^a_{zz})^{-1} \bm{\Sigma}^a_{zx}]^{-1}
  \]
  where
  \[
  \tilde{\bm{\Omega}}_n = \frac{1}{n_b}
  (\bm{\Sigma}^b_{zz})^{-1}\mathrm{Var}(y_i^b|\mathbf{z}_i^b) \ge \bm{\Omega}_n.
  \]
  Compared to \eqref{eq:tsiv-var}, it is larger
  than the variance of the efficient TSIV estimator. However, since
  the asymptotic variance of TSTSLS is larger than the efficient TSIV
  estimator, $\hat{\sigma}^2_{\text{naive}}(\hat{\beta}_{\text{TSTSLS}})$ may or may
  not over-estimate the variance of $\hat{\beta}_{\text{TSLS}}$. The naive
  variance estimator is used by \citet{gamazon2015gene} for gene
  testing. This is okay because under the null hypothesis $\beta_0 = 0$,
  we have $\tilde{\bm \Omega}_n = \bm \Omega_n$. However, the variance
  estimator is likely too small when constructing confidence intervals
  of $\beta$.
\end{remark}

\begin{remark}
  When $q = 1$, the covariance matrices all become scalars. The GMM estimator
  $\hat{\beta}_{n,\mathbf{W}}$ no longer depends on $\mathbf{W}$ and
  is always equal to the two-sample Wald ratio estimator. To see this,
  all the matrices in \eqref{eq:beta-hat} become scalars and
  \[
  \hat{\beta}_n = (s_{zy}^b/s_{zz}^b) / (s_{zx}^a/s_{zz}^a).
  \]
  The asymptotic variance of $\hat{\beta}_n$ is given by
  \eqref{eq:beta-hat-variance}, which can be simplified to
  \[
  \mathrm{Var}(\hat{\beta}_n) \approx \omega_n / (s_{zx}^a/s_{zz}^a)^2.
  \]
  The asymptotic variance in this special case can be derived more
  directly by the delta method as well \citep{burgess2015review}.
\end{remark}

\begin{remark} \label{rem:cov-mat}
  When $q > 1$, our results mean that the covariance matrices of
  $\mathbf{Z}$ are needed to compute any IV estimator and its
  asymptotic variance (unless only a single IV is used). Just
  observing the marginal regression coefficients
  is not enough. In situations where only the $\mathbf{S}_{zx}^a$ and
  $\mathbf{S}_{zy}^b$ are available (for example many GWAS only report
  summary statistics), one may estimate $\mathbf{S}_{zz}^a$ and
  $\mathbf{S}_{zz}^b$ (which reflects linkage disequilibrium in mendelian
  randomization) from additional datasets drawn from the same
  population. A similar idea of estimating linkage disequilibrium from
  additional dataset can be found in the context of multiple-SNP
  analysis in GWAS \citep{yang2012conditional}. In the context of
  Mendelian randomization, this means we can still compute the TSTSLS
  estimator by plugging in estimates of $\mathbf{\Sigma}_{zz}^a$ and
  $\mathbf{\Sigma}_{zz}^b$ obtained from other samples, but to compute
  the asymptotic variance, the matrix $\bm{\Omega}$ is not directly
  estimable because $\mathrm{Var}(y^b_i|\mathbf{z}^b_i)$ and $\mathrm{Var}(x^b_i|\mathbf{z}^b_i)$ are
  unknown. Nonetheless, one can still obtain a conservative estimate
  of $\bm{\Omega}$ from \eqref{eq:omega-2} using
  $\mathrm{Var}(y^b_i|\mathbf{z}^b_i) \le \mathrm{Var}(y^b_i)$ and
  $\mathrm{Var}(x^b_i|\mathbf{z}^b_i) \le \mathrm{Var}(x^b_i)$.
  This upper bound is usually not too conservative in Mendelian
  randomization since genetic variants identified so far usually only
  explain a small portion of the variability of complex diseases and
  traits \citep{manolio2009finding}.
\end{remark}


\section{Relaxing invariance and homogeneity assumptions}
\label{sec:relax-invar-assumpt}

Apart from the structural model and validity of IV
(\Cref{assump:npsem,assump:valid-iv}) that are necessary in the
one-sample setting, in \Cref{sec:tsiv-estim-using} we used additional
invariance/homogeneity and linearity assumptions
(\Cref{assump:invariance,assump:invariance-noise,assump:linearity}) to
identify and estimate the causal effect in the heterogeneous TSIV
setting. Next we attempt to relax these assumptions. Our main new
identification results in the next three sections are summarized in
\Cref{tab:summary}.

\renewcommand{\arraystretch}{1.25}
\begin{table}[t]
  \centering
  \caption{Summary of some identification results and assumptions made in
    this paper.}
  \label{tab:summary}
  \begin{tabularx}{\textwidth}{|l|X|c|c|c|}
    \hline
    \multirow{2}{*}{Assumption} & \multirow{2}{*}{Detail} & Prop.\ \ref{prop:1} & Prop.\ \ref{prop:2} &
    Prop.\
    \ref{prop:3}
    \\
    & & (Sec.\ \ref{sec:relax-invar-assumpt}) &
    (Sec.\
    \ref{sec:relax-line-instr})
    &
    (Sec.\ \ref{sec:relax-line-expos})\\
    \hline
    (\ref{assump:npsem}) Structural equation model & $y_i^s =
    g^s(x_i^s,u_i^s)$, $x_i^s = f^s(\mathbf{z}_i^s,v_i^s)$ &
    \checkmark & \checkmark & \checkmark
    \\
    (\ref{assump:valid-iv}) Validity of IV & $\mathbf{z}_i^s
    \independent (u_i^s,v_i^s)$ & \checkmark & \checkmark & \checkmark
    \\
    (\ref{assump:linearity}-1) Linearity of outcome eq.\ &
    $g^b(x_i,u_i) = \beta^b x_i + u_i$ & \checkmark & \checkmark & \\
    (\ref{assump:linearity}-2) Linearity of exposure eq.\ &
    $f^s(\mathbf{z}_i,v_i) = (\bm{\gamma}^s)^T \mathbf{z}_i + v_i$ &
    \checkmark & & \\
    (\ref{assump:invariance}) Structural invariance &
    $f^a = f^b$ & \checkmark & \checkmark & \checkmark  \\
    (\ref{assump:invariance-noise}) Sampling homogeneity of noise &
    $v_i^a \overset{d}{=} v^b_i$ &  &  & \checkmark  \\
    (\ref{assump:exposure-additive}) Additivity of exposure eq.\
    & $f^s(\mathbf{z},v) = f^s_z(\mathbf{z}) + f^s_v(v)$ & & \checkmark & \\
    (\ref{assump:monotonicity}) Monotonicity & $f^s(z,v)$ is
    monotone in $z$ & & & \checkmark \\
    \hline
    Identifiable estimand & & $\beta^b$ & $\beta^b$ & $\beta^b_{\mathrm{LATE}}$ \\
    \hline
  \end{tabularx}
\end{table}

First of all, notice that we did not use invariance of $g$ and $u$ in
the calculation above. Because $y^a$ is not observed, we do not need
to consider the exposure-outcome relation in sample $a$. In fact,
$u^a$ never appears in the calculation above, so we
can replace $\beta$ by $\beta^b$ and all the arguments in
\Cref{sec:tsiv-estim-using} still go through under the same
assumptions. For example, it is easy
to verify using \eqref{eq:gn-mean} that $\mathbf{m}_n(\beta^b)$ still has mean
$\mathbf{0}$ and converges to $\mathbf{0}$ in probability. Therefore,
the estimand of the TSIV estimators is indeed $\beta^b$ and we do not
need to assume $\beta^a = \beta^b$ or $u^a \overset{d}{=} u^b$. In
fact, $\beta^a$ is not identifiable from the data unless we link it to $\beta^b$.

Second, sampling homogeneity of the noise variable $v$
(\Cref{assump:invariance-noise}) is not crucially important in the
above linear structural equation models \eqref{eq:model-simple}. When
the expected values of $v^a$ and $v^b$ are different, they can be
absorbed in an intercept term and this does not affect the identification
and estimation of $\beta^b$. Also, our calculations above have already
considered the possibility that the variance of $v^a$ and $v^b$ are
different. To summarize, we have just shown that

\begin{proposition} \label{prop:1}
  Under Assumptions \ref{assump:npsem}, \ref{assump:valid-iv},
  \ref{assump:linearity}-1 (for sample $b$),
  \ref{assump:linearity}-2 (for both samples) and
  \ref{assump:invariance}, the TSIV estimators in
  \Cref{sec:tsiv-estim-using} can consistently estimate $\beta^b$.
\end{proposition}

Thus noise homogeneity (\Cref{assump:invariance-noise}) is not necessary
when the structural relations are linear. However, we will see in the
next two Sections that \Cref{assump:invariance-noise} is quite
important when the structural relations are not linear.


\section{Relaxing linearity of the instrument-exposure equation}
\label{sec:relax-line-instr}

In one-sample IV analyses, correct specification of
the instrument-exposure model is not necessary for consistent
estimation of the causal effect. To see this, suppose the linear
exposure-outcome model is correctly specified in
\eqref{eq:model-simple} (i.e.\ Assumption \ref{assump:linearity}-2
holds). In the one-sample problem, the parameter $\beta$ can be
identified by the following moment condition
\[
\mathbb{E}[h(\mathbf{z}) \cdot (y - x \beta)] = 0
\]
for any function $h(\mathbf{z})$ due to the independence of
$\mathbf{z}$ and $u$ as long as $\mathrm{Cov}(x, h(\mathbf{z})) \ne
0$. This results in the class of instrumental
variable estimators
\begin{equation} \label{eq:iv-estimator}
  \hat{\beta}_h = \Big[\sum_{i=1}^n y_i
  h(\mathbf{z}_i)\Big] \Big/ \Big[\sum_{i=1}^n x_i h(\mathbf{z}_i)\Big].
\end{equation}
The TSLS estimator is a special case of \eqref{eq:iv-estimator} when
$h(\mathbf{z}) = \mathbf{z}^T \bm{\gamma}$ and $\bm{\gamma}$ is estimated from the first
stage regression. In general, $\hat{\beta}_h$ is consistent and
asymptotically normal. The asymptotic variance of $\hat{\beta}_h$ depends on the choice
of $h$. The optimal choice of $h$, often called the optimal
instrument, is the conditional expectation of $x$ given $\mathbf{z}$:
$h^{*}(\mathbf{z}) = \mathrm{E}[x|\mathbf{z}]$. To summarize, in the
one-sample problem, the TSLS estimator is consistent for $\beta$ even
if the linear instrument-exposure model is misspecified, although in
that case the TSLS estimator may be less efficient than the optimal instrumental
variable estimator. We refer the reader to
\citet{vansteelandt2015robustness} for a recent discussion on
robustness and efficiency of one-sample IV estimators under various
types of model misspecification.

\begin{figure}[t]
  \centering
  \begin{subfigure}[t]{0.48\textwidth}
    \includegraphics[width = \textwidth]{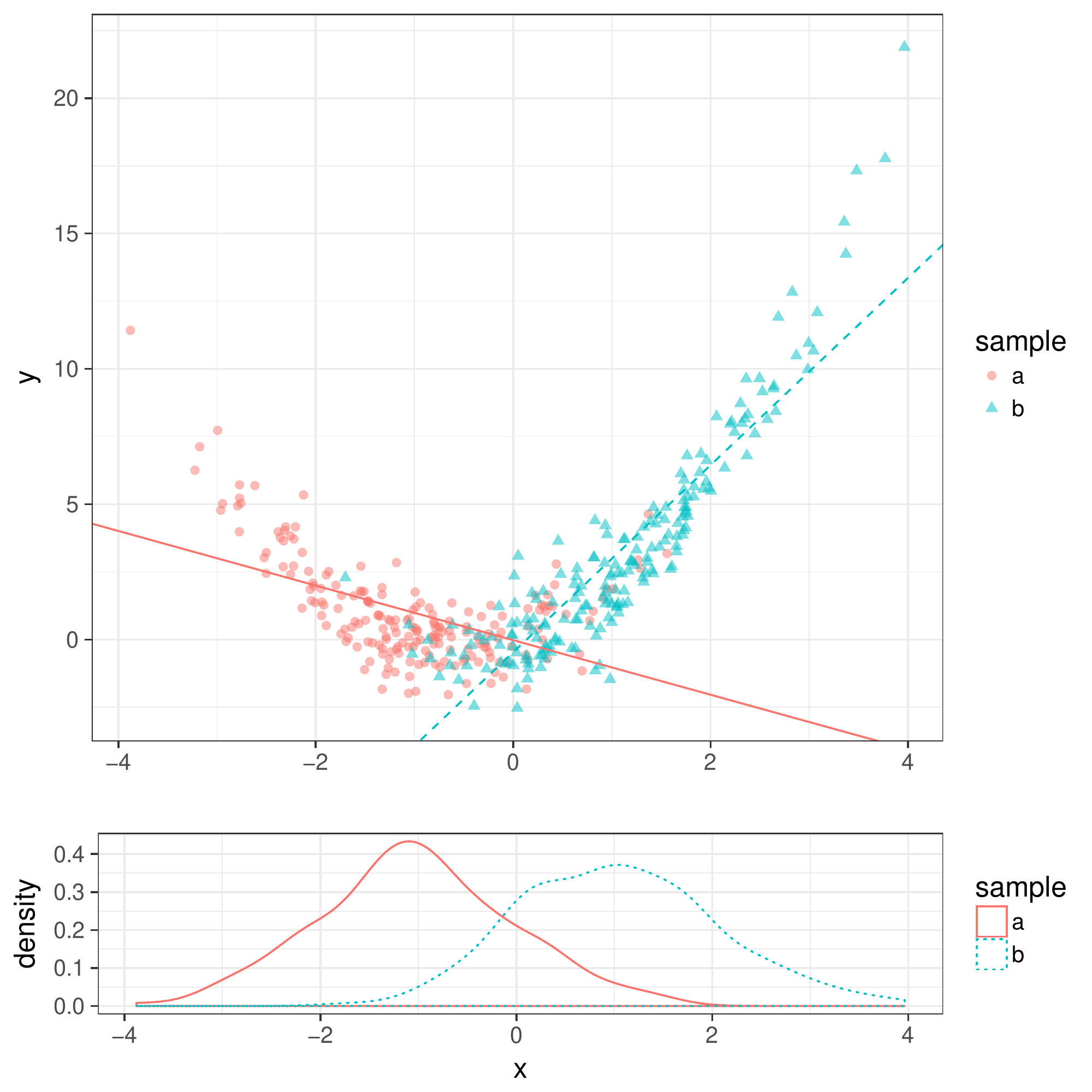}
    \caption{Before matching.}
  \end{subfigure}
  \begin{subfigure}[t]{0.48\textwidth}
    \includegraphics[width = \textwidth]{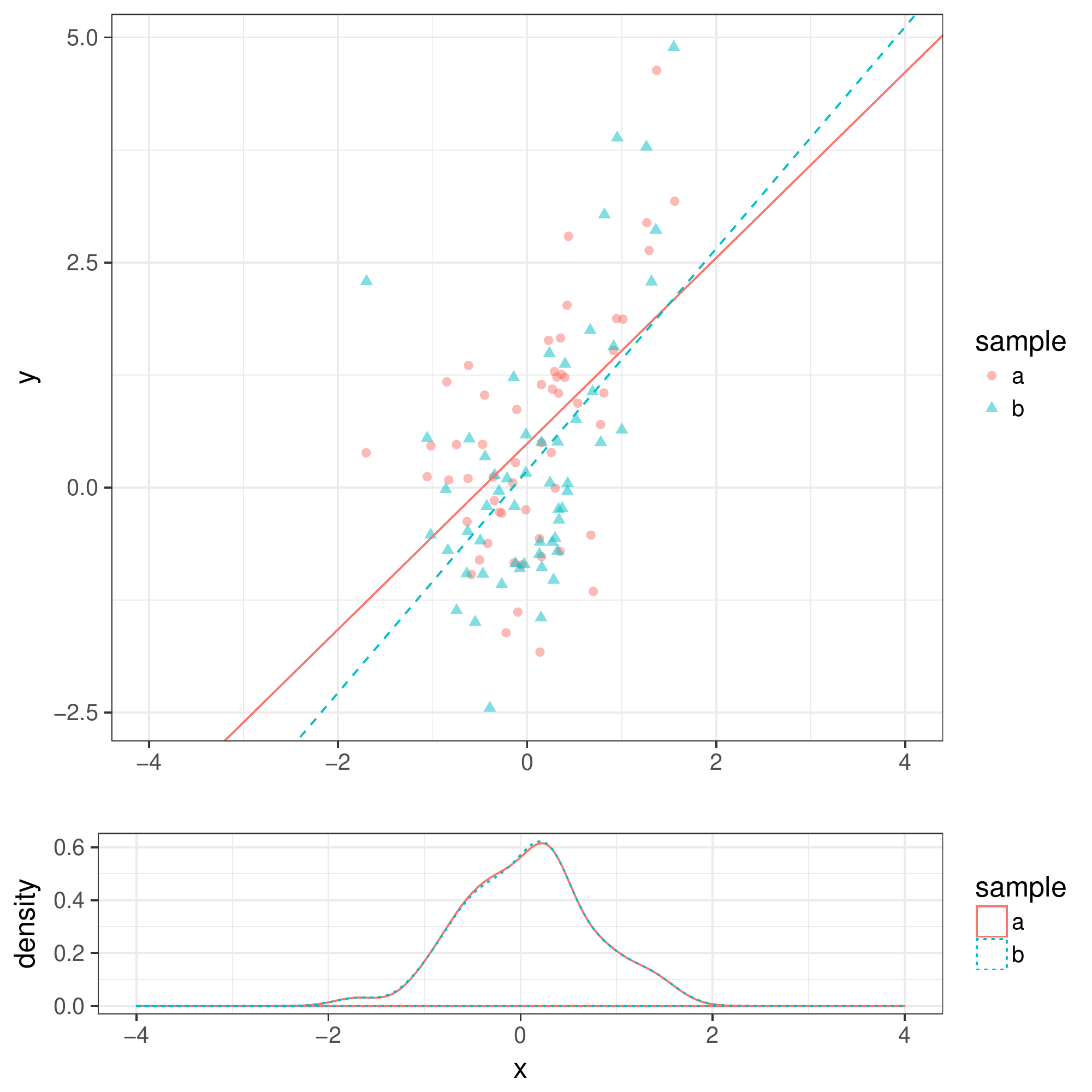}
    \caption{After nearest neighbor matching (pairs are kept if the
      absolute within pair difference is less than $0.1$).}
  \end{subfigure}
  \caption{``Conspiracy'' of model misspecification and random
    design. In this example, $f^a(z,v) = f^b(z,v) = z^2 + z + v$ where
    $v^a \overset{d}{=} v^b \sim \mathrm{N}(0, 1)$, $z^a \sim
    \mathrm{N}(-1,1)$, and $z^b \sim \mathrm{N}(1,1)$. If the model is
    misspecified and a linear model is used as an approximation,
    $f^s(z, v) \approx \eta^s + \gamma^s z^s$, the projections
    $\gamma^a$ and $\gamma^b$ depends on the distribution of $z^a$ and
    $z^b$ and have different signs in this example. By only considering
    the common support of the two samples and matching the observations,
    the projections $\gamma^a$ and $\gamma^b$ are much closer.}
  \label{fig:conspiracy}
\end{figure}

This robustness property of TSLS does not carry to the two-sample
setting due to a phenomenon known as the ``conspiracy'' of model
misspecification and random design \citep{white1980using,buja2014models}. Under the
general instrument-exposure equation $x_i^s = f^s(\mathbf{z}_i^s,v_i^s)$
in \eqref{eq:npsem-exposure}, the best linear projection (in Euclidean distance)
\begin{equation} \label{eq:gamma-projection}
  \bm{\gamma}^s = \argmin_{\bm{\gamma}} \mathrm{E}\big\{[(\mathbf{z}_i^s)^T
  \bm{\gamma} - f^s(\mathbf{z}_i^s,v_i^s)]^2\big\}
\end{equation}
depends on the structural function $f^s$, the distribution of the noise
variable $v^s$, and the distribution of the instrumental variables
$\mathbf{z}_i^s$. Therefore, even if structural invariance
(\Cref{assump:invariance}) and sampling homogeneity of the noise variables
(\Cref{assump:invariance-noise}) are satisfied, the best linear
approximations $\bm{\gamma}^a$ and
$\bm{\gamma}^b$ can still be different if the sampling distributions
of $\mathbf{z}$ are different. In extreme cases, $\bm{\gamma}^a$ and
$\bm{\gamma}^b$ can even have different signs; see
\Cref{fig:conspiracy} for an example. Since the TSTSLS estimator
converges to $\gamma^b \beta^b / \gamma^a$ when the instrumental
variable is univariate, the TSTSLS estimator and other TSIV estimators
are biased and may even estimate the sign of $\beta^b$ incorrectly.

There are two ways to mitigate the issue of non-linearity of the
instrument-exposure equation. The first is to only consider the common
support of $\mathbf{z}^a$ and $\mathbf{z}^b$ as suggested by
\citet{lawlor2016commentary} and match or weight the observations so that
$\mathbf{z}^a$ and $\mathbf{z}^b$ have the same distribution. This
ensures the projections $\bm{\gamma}^a$ and $\bm{\gamma}^b$ are the
same and is illustrated in \Cref{fig:conspiracy}. The second solution is to nonparametrically model the
instrument-exposure relation to avoid the drawback of using the linear
approximations. However, this is difficult if the dimension of the IVs
is high.

We want to emphasize that, unlike the scenario with linear
instrument-exposure equation in \Cref{sec:tsiv-estim-using}, both
solutions above still hinge on sampling
invariance of noise variables (\Cref{assump:invariance-noise}). Even
if the distributions of $\mathbf{z}^a$ and $\mathbf{z}^b$ are the same
and $f^a$ is modeled nonparametrically, the best linear or nonlinear
approximation still depends on the distribution of the noise variable
$v$. If \Cref{assump:invariance-noise} is violated so $v^a$ and $v^b$
have different distributions, the TSIV estimators are still generally
biased, though the bias is unlikely to be extremely large. It is also worth
noting that sampling homogeneity of the noise variables
(\Cref{assump:invariance-noise}) is untestable in the two sample
setting because $\mathbf{x}^b$ is not observed.

One way to relax \Cref{assump:invariance-noise} is to assume the
instrument-exposure equation is additive:
\begin{assumption} \label{assump:exposure-additive}
  (Additivity of the instrument-exposure equation) $f^s(\mathbf{z},v)
  = f^s_z(\mathbf{z}) + f^s_v(v)$.
\end{assumption}
Under \Cref{assump:exposure-additive}, we may non-parametrically
estimate $f_z^s(\mathbf{z})$ and then estimate $\beta^b$ by
regressing $y_i^b$ on the predicted $f_z^a(\mathbf{z}_i^b)$. This is
consistent for $\beta^b$ if $f_z^a$ is estimated consistently, because
\[
y_i^b = \beta^b x_i^b + u_i^b = \beta^b f^b_z(\mathbf{z}_i^b) + \beta^b
f^b_v(v_i^b) + u_i^b = \beta^b f^a_z(\mathbf{z}_i^b) + (\beta^b
f^a_v(v_i^b) + u_i^b).
\]
The last equation used structural invariance
(\Cref{assump:invariance}). Even if the noise variables $u$ and $v$
may have different distributions in the two samples, the estimation of
$\beta^b$ is not affected (see \Cref{sec:relax-invar-assumpt}). To
summarize, we have shown that
\begin{proposition} \label{prop:2}
  In \Cref{prop:1} and absence of noise homogeneity, $\beta^b$ can still
  be identified if the exposure equation is additive.
\end{proposition}

\section{Relaxing linearity of the exposure-outcome equation}
\label{sec:relax-line-expos}

\subsection{LATE in the one-sample setting}
\label{sec:brief-review-one}

When the exposure-outcome equation is nonlinear, an additional
assumption called homogeneity is usually needed to identify the causal
effect. Next we review this approach in the one-sample setting when
the instrument and the exposure are both binary. In this case, we can
define four classes of
observations based on the instrument-exposure equation: for $s = a,b$,
\[
t^s(v) =
\begin{cases}
  \mathrm{always~taker~(at)}, & \mathrm{if}~f^s(0,v) = 1,~f^s(1,v)=1,
  \\
  \mathrm{complier~(co)}, & \mathrm{if}~f^s(0,v) = 0,~f^s(1,v)=1, \\
  \mathrm{never~taker~(nt)}, & \mathrm{if}~f^s(0,v) = 0,~f^s(1,v)=0, \\
  \mathrm{defier~(de)}, & \mathrm{if}~f^s(0,v) = 1,~f^s(1,v)=0. \\
\end{cases}
\]
Classes are important to remove endogeneity since conditioning on the
class, the exposure $x$ is no longer dependent on the noise variable
$u$, that is
\begin{equation}
  \label{eq:class-independence}
  x^s \independent u^s \, | \, t^s(v^s).
\end{equation}
The last equation is true because given $t^s(v^s)$ and hence
the values of $f^s(0,v)$ and $f^s(1,v)$, the only randomness of $x^s$
comes from $z^s$ which is independent of $u^s$. If the classes were
observable, \eqref{eq:class-independence} implies that we can identify
the class-conditional average outcome $\mathrm{E}[g^s(x,u^s)|t^s=t]$
for $(t,x)$ in the support of $(t^s,x^s)$,
which is a subset of $\{at,co,nt,de\} \times \{0, 1\}$. More
specifically, since $\mathrm{P}(x^s=0|t^s=at) = 0$ and
$\mathrm{P}(x^s=1|t^s=nt) = 0$, the support of $(t^s,x^s)$ contains
$6$ elements,
$\mathrm{supp}(t^s,x^s)=\{(at,1),(co,0),(co,1),(nt,0),(de,0),(de,1)\}$. However,
the classes are not directly observable, and in fact we can only identify
four conditional expectations $\mathrm{E}[y^s|x^s=x,z^s=z] =
\mathrm{E}[g^s(x,u^s)|x^s=x,z^s=z]$ from the data. This means that the
class-conditional average outcomes are not identifiable, because in the
following system of equations,
\begin{equation} \label{eq:class-ate-identification}
  \begin{split}
    \mathrm{E}[g^s(0,u^s)|x^s=0,z^s=0] &= \mathrm{E}[g^s(0,u^s)|t^s = nt] \cdot
    \mathrm{P}(t^s = nt) + \mathrm{E}[g^s(0,u^s)|t^s = co] \cdot
    \mathrm{P}(t^s = co), \\
    \mathrm{E}[g^s(0,u^s)|x^s=0,z^s=1] &= \mathrm{E}[g^s(0,u^s)|t^s = nt] \cdot
    \mathrm{P}(t^s = nt) + \mathrm{E}[g^s(0,u^s)|t^s = de] \cdot
    \mathrm{P}(t^s = de), \\
    \mathrm{E}[g^s(1,u^s)|x^s=1,z^s=0] &= \mathrm{E}[g^s(1,u^s)|t^s = at] \cdot
    \mathrm{P}(t^s = at) + \mathrm{E}[g^s(1,u^s)|t^s = de] \cdot
    \mathrm{P}(t^s = de), \\
    \mathrm{E}[g^s(1,u^s)|x^s=1,z^s=1] &= \mathrm{E}[g^s(1,u^s)|t^s = at] \cdot
    \mathrm{P}(t^s = at) + \mathrm{E}[g^s(1,u^s)|t^s = co] \cdot
    \mathrm{P}(t^s = co),
  \end{split}
\end{equation}
there are six class-conditional average outcomes but only four
equations. Note that to derive \eqref{eq:class-ate-identification} we have
used \Cref{assump:valid-iv} which asserts $z^s \independent t^s = t^s(v^s)$ and
$z^s \independent u^s$, so $\mathrm{E}[g^s(x,u^s)|z^s,t^s] =
\mathrm{E}[g^s(x,u^s)|t^s]$ and $\mathrm{P}(t^s = t|z^s) =
\mathrm{P}(t^s = t)$ for any fixed $x$ and $t$.

The monotonicity assumption is used to reduce the number of free
parameters in \eqref{eq:class-ate-identification}.
\begin{assumption} \label{assump:monotonicity}
  (Monotonicity) $f^s(z,v)$ is a monotone function of $z$ for any $v$ and $s =
  a,b$.
\end{assumption}
Without loss of generality, we will assume $f^s(z,v)$ is an increasing
function of $z$, otherwise we can use $-x^s = -f^s(z^s,v^s)$ as the
exposure. In the context of binary instrument and binary exposure,
\Cref{assump:monotonicity} means that $\mathrm{P}(t^s = de) = 0$ and is
often called the \emph{no-defiance} assumption
\citep{balke1997bounds}. This eliminates two class-conditional average
outcomes, $\mathrm{E}[g^s(0,u^s)|t^s=de]$ and
$\mathrm{E}[g^s(1,u^s)|t^s=de]$, leaving us four equations and four
class-conditional average outcomes. Therefore, using
\eqref{eq:class-ate-identification}, we can identify the so called
\emph{local average treatment effect} (LATE),
$\mathrm{E}[g^s(1,u^s)-g^s(0,u^s)|t^s=co]$
\citep{angrist1996identification}. In particular, under
\Cref{assump:npsem,assump:valid-iv,assump:monotonicity}, one can show
that the TSLS estimator in sample $s$ converges to
\begin{equation} \label{eq:tsls-estimand}
  \begin{split}
    \beta_{\mathrm{LATE}}^s=&\frac{\mathrm{E}[y^s|z^s=1] - \mathrm{E}[y^s|z^s =
      0]}{\mathrm{E}[x^s|z^s=1] - \mathrm{E}[x^s|z^s = 0]} \\
    =&
    \frac{\mathrm{E}[g^s(1,u^s)-g^s(0,u^s)|t^s=co] \cdot \mathrm{P}(t^s =
      co)}{\mathrm{P}(t^s = co)} \\
    =&
    \mathrm{E}[g^s(1,u^s)-g^s(0,u^s)|t^s=co].
  \end{split}
\end{equation}
See \eqref{eq:tstsls-estimand} below for proof this result.

When the exposure $x$ is continuous, we may still define the class $t$
such that \eqref{eq:class-independence} holds and identify the
class-conditional average outcomes on the joint support of $x$ and
$t$. This support may be very limited when the instrument $z$ is
binary. We refer the reader to \citet{imbens2007nonadditive} for
further detail and discussion. In this case, the instrumental variable
estimator $\hat{\beta}_h$ in \eqref{eq:iv-estimator} converges in
probability to a weighted average of local
average treatment effects \citep{angrist2000interpretation}. Note that
in order for the weights to be non-negative, ordering the instruments
by $\mathrm{E}[x^s|z^s=z]$ must simultaneously order the instruments by
the value of $h(z)$ \citep[Theorem 2,3]{angrist2000interpretation}. A
preferable choice of $h(z)$ is the conditional expectation
$\mathrm{E}[x^s|z^s=z]$.

\subsection{LATE in the two-sample setting}
\label{sec:late-two-sample}

We can still follow the LATE framework in the two-sample setting
considered in this paper. When the instrument and the exposure
are both binary, the TSTSLS estimator converges to a modification of
\eqref{eq:tsls-estimand} by taking the expectations in the numerator
over sample $a$ and the expectations in the denominator over sample
$b$,
\begin{equation}
  \label{eq:tstsls-estimand}
  \begin{split}
    \beta_{\mathrm{LATE}}^{ab}=&\frac{\mathrm{E}[y^b|z^b=1] - \mathrm{E}[y^b|z^b =
      0]}{\mathrm{E}[x^a|z^a=1] - \mathrm{E}[x^a|z^a = 0]} \\
    =&
    \frac{\mathrm{E}[g^b(1,u^b)-g^b(0,u^b)|t^b=co] \cdot \mathrm{P}(t^b =
      co)}{\mathrm{P}(t^a = co)} \\
    =& \beta_{\mathrm{LATE}}^b \cdot \frac{\mathrm{P}(t^b =
      co)}{\mathrm{P}(t^a = co)}.
  \end{split}
\end{equation}
Next we prove
the second equality in \eqref{eq:tstsls-estimand}. First we consider the
numerator
\[
\begin{split}
  &\mathrm{E}[y^b|z^b=1] - \mathrm{E}[y^b|z^b = 0] \\
  =&\sum_{t \in \{at,co,nt,de\}} \big(\mathrm{E}[y^b|z^b=1,t^b=t] -
  \mathrm{E}[y^b|z^b = 0,t^b=t]\big) \cdot \mathrm{P}(t^b=t) \\
  =&\sum_{t \in \{at,co,nt\}} \big(\mathrm{E}[y^b|z^b=1,t^b=t] -
  \mathrm{E}[y^b|z^b = 0,t^b=t]\big) \cdot \mathrm{P}(t^b=t) \\
\end{split}
\]
where the first equality is due to the law of total expectation and
the second equality uses \Cref{assump:monotonicity}. Next, notice that
$y^b \independent z^b |t^b = at$, because $\mathrm{P}(x^b=1|t^b=at) =
1$ and by the exclusion restriction (implied from
\Cref{assump:valid-iv}), $y^b$ only depends on $z^b$ through
$x^b$. Similarly, $y^b \independent z^b |t^b = ne$. Therefore, we are
left with just the compliers
\[
\begin{split}
  &\mathrm{E}[y^b|z^b=1] - \mathrm{E}[y^b|z^b = 0] \\
  =&\big(\mathrm{E}[y^b|z^b=1,t^b=co] -
  \mathrm{E}[y^b|z^b = 0,t^b=co]\big) \cdot \mathrm{P}(t^b=co) \\
  =&\big(\mathrm{E}[g^b(1,u^b)|z^b=1,t^b=co] -
  \mathrm{E}[g^b(0,u^b)|z^b = 0,t^b=co]\big) \cdot
  \mathrm{P}(t^b=co) \\
  =&\big(\mathrm{E}[g^b(1,u^b)|t^b=co] -
  \mathrm{E}[g^b(0,u^b)|t^b=co]\big) \cdot \mathrm{P}(t^b=co). \\
\end{split}
\]
In the last equation we have again used the exclusion
restriction. Similarly, the denominator in \eqref{eq:tstsls-estimand}
is
\[
\begin{split}
  &\mathrm{E}[x^b|z^b=1] - \mathrm{E}[x^b|z^b = 0] \\
  =&\sum_{t \in \{at,co,nt\}} \big(\mathrm{E}[x^b|z^b=1,t^b=t] -
  \mathrm{E}[x^b|z^b = 0,t^b=t]\big) \cdot \mathrm{P}(t^b=t) \\
  =& (1 - 1) \cdot \mathrm{P}(t^b = at) + (1 - 0) \cdot
  \mathrm{P}(t^b=co) + (0 - 0) \cdot \mathrm{P}(t^b=nt) \\
  =& \mathrm{P}(t^b = co).
\end{split}
\]
Finally, note that similar to the one-sample case,
\eqref{eq:tstsls-estimand} only uses
\Cref{assump:npsem,assump:valid-iv,assump:monotonicity}.

When
structural invariance of the instrument-exposure equation $f$
(\Cref{assump:invariance}) and sampling homogeneity of the noise variable $v$
(\Cref{assump:invariance-noise}) hold, we have $t^a
\overset{d}{=} t^b$ because the class $t^s$ is a function of
$f^s$ and $v^s$. Therefore $\beta_{\mathrm{LATE}}^{ab} = \beta_{\mathrm{LATE}}^b$
by equation \eqref{eq:tstsls-estimand}. To summarize, we have just
shown that
\begin{proposition} \label{prop:3}
  When there is one binary instrument and one binary exposure, under
  Assumptions \ref{assump:npsem}, \ref{assump:valid-iv},
  \ref{assump:invariance-noise} and \ref{assump:monotonicity}, the TSTSLS estimator identifies $\beta^b_{\text{LATE}}$.
\end{proposition}

In general, the estimand of TSTSLS is a
scaling of the LATE in the sample $b$. Since $f^a$ and $f^b$ are
non-trivial functions of $z$ by \Cref{assump:valid-iv}, the
proportions of compliers are positive and hence the ratio $\mathrm{P}(t^b =
co)/\mathrm{P}(t^a = co) > 0$. This means that
$\beta_{\mathrm{LATE}}^{ab}$ has the same sign as
$\beta_{\mathrm{LATE}}^b$.

When the exposure is continuous, most of the arguments in
\citet{angrist2000interpretation} would still hold as they were proved
separately for the numerator and the denominator just like our proof of
\eqref{eq:tstsls-estimand}. Similarly, the TSTSLS estimator
converges in probability to the estimand of the TSLS estimator in
sample $b$ times a scaling factor, and the scaling factor is equal to
$1$ under \Cref{assump:invariance,assump:invariance-noise}. However,
the scaling factor is not always positive because in the absence of
\Cref{assump:invariance-noise}, the
conditional expectation $\mathrm{E}[x^s|z^s]$ can be different in the
two samples (same issue as in \Cref{sec:relax-line-instr}). Similar to
\Cref{sec:relax-line-instr}, this can be resolved by assuming
additivity (\Cref{assump:exposure-additive}).

\section{Simulation}
\label{sec:simulation}

We evaluate the efficiency and robustness of the linear TSIV
estimators using numerical simulation. In all simulations we consider
$10$ binary instrumental variables generated by
\begin{equation} \label{eq:sim-iv}
  \mathbf{z}_i^s = \mathrm{sign}(\mathbf{z}_i^{s*}),~\mathbf{z}_i^{s*}
  \overset{i.i.d.}{\sim} \mathrm{N}(\mathbf{1},
  \bm{\Sigma}^s),~\Sigma_{jk}^s = (\rho^s)^{|j-k|},~s=a,b.
\end{equation}

We first verify the asymptotic results regarding the TSIV estimators in
\Cref{sec:tsiv-estim-using}. In our first simulation, the exposures
and the outcomes
are generated by
\begin{align}
  x_i^s &= 0.2 \cdot (\mathbf{1}^T \mathbf{z}_i^s) +
  v_i^s, \label{eq:sim1-exposure} \\
  y_i^s &= x_i^s + u_i^s, \label{eq:sim1-outcome}\\
  (v_i^s,u_i^s) &\overset{i.i.d.}{\sim} \mathrm{N}\left(
    \begin{pmatrix}
      0 \\
      0 \\
    \end{pmatrix},
    \begin{pmatrix}
      1 & \sigma_{uv} \\
      \sigma_{uv} & 1 \\
    \end{pmatrix}
  \right),~i=1,\dotsc,n^s,~s=a,b. \label{eq:sim1-uv}
\end{align}
In this simulation we used $\rho^a = 0.5$, $\rho^b = 0.5$, $0$, or
$-0.5$, $n^a=1000$ or $5000$, $n^b = 1000$ or $5000$, and $\sigma_{uv}
= 0.5$.

In \Cref{tab:sim1}, we compare the performance of the TSTSLS
estimator and the optimal TSIV estimator after centering the
variables. In particular, we report
the bias, standard deviation (SD), average standard error (SE), and coverage of the 95\%
asymptotic confidence interval. When $\rho^a = \rho^b = 0.5$, the two
estimators are asymptotically equivalent by \eqref{eq:omega-2}. This
is verified by \Cref{tab:sim1} as the two estimators have the same
bias, variance, and coverage in this case. When $\rho^a$ and $\rho^b$ are
different, the optimal TSIV estimator should be more efficient
than TSTSLS (at least theoretically). In the simulations we find that
in almost all cases the two estimators have the same variance, but the
optimal TSIV estimator has smaller finite sample bias.
The difference between the optimal TSIV estimator and the TSTSLS
estimator is
substantial only if $\bm{\Sigma}^a$ and $\bm{\Sigma}^b$ (in this
simulation, $\rho^a$ and $\rho^b$) are very different and $n^b$ is
much larger than $n^a$. This phenomenon can also be seen from
\eqref{eq:omega-2} as discussed in \Cref{rem:efficiency}.

\renewcommand{\arraystretch}{1.25}
\begin{table}[t]
  \centering
  \caption{Simulation 1: Asymptotic efficiency of the TSTSLS and optimal TSIV
    estimators. The reported numbers are obtained by 10000 realizations
    of data simulated from equations \eqref{eq:sim-iv},
    \eqref{eq:sim2-exposure}, \eqref{eq:sim1-outcome}, and
    \eqref{eq:sim1-uv}.}
  \label{tab:sim1}
  \begin{tabular}{|lllll|rrrr|rrrr|}
    \hline
    & & & & & \multicolumn{4}{c|}{TSTSLS} &
    \multicolumn{4}{c|}{Optimal TSIV} \\
    $\beta$ & $\rho^a$ & $\rho^b$ & $n^a$ & $n^b$ & Bias & SD & SE & Cover &
    Bias & SD & SE & Cover \\
    \hline
    1 & 0.5 & 0.5 & 1000 & 1000  & $-0.020$ & $0.100$ & $0.100$ & $0.941$ & $-0.020$ & $0.100$ & $0.100$ & $0.941$ \\
    &  &  &  & 5000  & $-0.021$ & $0.062$ & $0.063$ & $0.924$ & $-0.021$ & $0.062$ & $0.063$ & $0.923$ \\
    &  &  & 5000 & 1000  & $-0.006$ & $0.090$ & $0.090$ & $0.949$ & $-0.006$ & $0.090$ & $0.090$ & $0.949$ \\
    &  &  &  & 5000  & $-0.004$ & $0.045$ & $0.045$ & $0.948$ & $-0.004$ & $0.045$ & $0.045$ & $0.948$ \\
    &  & 0 & 1000 & 1000  & $-0.046$ & $0.126$ & $0.126$ & $0.928$ & $-0.039$ & $0.127$ & $0.126$ & $0.933$ \\
    &  &  &  & 5000  & $-0.047$ & $0.072$ & $0.072$ & $0.875$ & $-0.029$ & $0.072$ & $0.071$ & $0.909$ \\
    &  &  & 5000 & 1000  & $-0.012$ & $0.121$ & $0.120$ & $0.948$ & $-0.011$ & $0.121$ & $0.120$ & $0.947$ \\
    &  &  &  & 5000  & $-0.010$ & $0.057$ & $0.057$ & $0.947$ & $-0.008$ & $0.057$ & $0.057$ & $0.948$ \\
    &  & -0.5 & 1000 & 1000  & $-0.058$ & $0.135$ & $0.135$ & $0.918$ & $-0.045$ & $0.137$ & $0.134$ & $0.924$ \\
    &  &  &  & 5000  & $-0.058$ & $0.077$ & $0.075$ & $0.851$ & $-0.031$ & $0.076$ & $0.074$ & $0.907$ \\
    &  &  & 5000 & 1000  & $-0.012$ & $0.130$ & $0.129$ & $0.951$ & $-0.011$ & $0.130$ & $0.129$ & $0.951$ \\
    &  &  &  & 5000  & $-0.013$ & $0.060$ & $0.061$ & $0.947$ & $-0.010$
    & $0.061$ & $0.061$ &
    $0.949$ \\
    \hline
    10 & 0.5 & 0.5 & 1000 & 1000  & $-0.197$ & $0.726$ & $0.726$ & $0.929$ & $-0.197$ & $0.724$ & $0.724$ & $0.929$ \\
    &  &  &  & 5000  & $-0.194$ & $0.548$ & $0.548$ & $0.915$ & $-0.192$ & $0.544$ & $0.545$ & $0.916$ \\
    &  &  & 5000 & 1000  & $-0.036$ & $0.577$ & $0.578$ & $0.949$ & $-0.036$ & $0.577$ & $0.578$ & $0.949$ \\
    &  &  &  & 5000  & $-0.037$ & $0.327$ & $0.327$ & $0.948$ & $-0.037$ & $0.327$ & $0.327$ & $0.948$ \\
    &  & 0 & 1000 & 1000  & $-0.468$ & $0.867$ & $0.866$ & $0.898$ & $-0.330$ & $0.870$ & $0.860$ & $0.920$ \\
    &  &  &  & 5000  & $-0.475$ & $0.596$ & $0.585$ & $0.836$ & $-0.249$ & $0.587$ & $0.574$ & $0.900$ \\
    &  &  & 5000 & 1000  & $-0.096$ & $0.755$ & $0.754$ & $0.947$ & $-0.086$ & $0.756$ & $0.753$ & $0.948$ \\
    &  &  &  & 5000  & $-0.102$ & $0.393$ & $0.393$ & $0.938$ & $-0.072$ & $0.394$ & $0.392$ & $0.941$ \\
    &  & -0.5 & 1000 & 1000  & $-0.586$ & $0.932$ & $0.915$ & $0.876$ & $-0.380$ & $0.933$ & $0.902$ & $0.907$ \\
    &  &  &  & 5000  & $-0.575$ & $0.626$ & $0.610$ & $0.808$ & $-0.254$ & $0.598$ & $0.585$ & $0.902$ \\
    &  &  & 5000 & 1000  & $-0.112$ & $0.807$ & $0.807$ & $0.948$ & $-0.093$ & $0.808$ & $0.807$ & $0.948$ \\
    &  &  &  & 5000  & $-0.118$ & $0.420$ & $0.415$ & $0.934$ & $-0.071$ & $0.419$ & $0.413$ & $0.943$ \\
    \hline
  \end{tabular}
\end{table}

In the second simulation, we examine how misspecification of the
instrument-exposure equation may bias the TSIV estimator. The data are
generated in the same way as in the first simulation except that we
add interaction terms in the instrument-exposure equation. More
specifically, \eqref{eq:sim1-exposure} is replaced by
\begin{equation}
  x_i^s = 0.2 \cdot (\mathbf{1}^T \mathbf{z}_i^s) + 0.02 \cdot
  \sum_{j \ne k} z_{ij}^s z_{ik}^s + v_i^s. \label{eq:sim2-exposure}
\end{equation}

The results of the second simulation are reported in
\Cref{tab:sim2}. When $\rho^a = \rho^b = 0.5$, the TSTSLS and the optimal
TSIV estimators are still unbiased and the confidence intervals
provide desired coverage. This is because the best linear
approximations of the instrument-exposure equation are the same in the
two samples. However, when $\rho^a \ne \rho^b$, the TSTSLS and the
optimal TSIV estimators are biased and failed to cover the true
parameter at the nominal $95\%$ rate. As discussed in
\Cref{sec:relax-line-instr}, this is because the best linear
approximations of the instrument-exposure equation are different in
the two samples. In addition, note that the optimal TSIV estimator
tends to have larger bias in this simulation.

\begin{table}[t]
  \caption{Simulation 2: When the instrument-exposure equation is misspecified, the
    TSIV estimators can be biased. The reported numbers are obtained by 10000 realizations
    of data simulated from equations \eqref{eq:sim-iv},
    \eqref{eq:sim3-exposure}, \eqref{eq:sim1-outcome}, and
    \eqref{eq:sim3-uv}.}
  \label{tab:sim2}
  \centering
  \begin{tabular}{|lllll|rrrr|rrrr|}
    \hline
    & & & & & \multicolumn{4}{c|}{TSTSLS} &
    \multicolumn{4}{c|}{Optimal TSIV} \\
    $\beta$ & $\rho^a$ & $\rho^b$ & $n^a$ & $n^b$ & Bias & SD & SE & Cover &
    Bias & SD & SE & Cover \\
    \hline
    1 & 0.5 & 0.5 & 1000 & 1000  & $-0.009$ & $0.069$ & $0.067$ & $0.938$ & $-0.009$ & $0.068$ & $0.067$ & $0.940$ \\
    &  &  &  & 5000  & $-0.011$ & $0.044$ & $0.042$ & $0.922$ & $-0.010$ & $0.044$ & $0.042$ & $0.922$ \\
    &  &  & 5000 & 1000  & $-0.001$ & $0.061$ & $0.060$ & $0.946$ & $-0.001$ & $0.061$ & $0.060$ & $0.946$ \\
    &  &  &  & 5000  & $-0.002$ & $0.031$ & $0.030$ & $0.942$ & $-0.002$ & $0.031$ & $0.030$ & $0.942$ \\
    &  & 0 & 1000 & 1000  & $\phantom{-}0.041$ & $0.086$ & $0.085$ & $0.927$ & $\phantom{-}0.046$ & $0.086$ & $0.085$ & $0.920$ \\
    &  &  &  & 5000  & $\phantom{-}0.042$ & $0.051$ & $0.050$ & $0.878$ & $\phantom{-}0.052$ & $0.051$ & $0.049$ & $0.827$ \\
    &  &  & 5000 & 1000  & $\phantom{-}0.059$ & $0.081$ & $0.080$ & $0.885$ & $\phantom{-}0.060$ & $0.081$ & $0.080$ & $0.884$ \\
    &  &  &  & 5000  & $\phantom{-}0.060$ & $0.039$ & $0.038$ & $0.665$ & $\phantom{-}0.061$ & $0.039$ & $0.038$ & $0.651$ \\
    &  & -0.5 & 1000 & 1000  & $\phantom{-}0.050$ & $0.092$ & $0.091$ & $0.919$ & $\phantom{-}0.059$ & $0.093$ & $0.091$ & $0.904$ \\
    &  &  &  & 5000  & $\phantom{-}0.051$ & $0.054$ & $0.052$ & $0.847$ & $\phantom{-}0.067$ & $0.054$ & $0.052$ & $0.755$ \\
    &  &  & 5000 & 1000  & $\phantom{-}0.074$ & $0.086$ & $0.086$ & $0.860$ & $\phantom{-}0.075$ & $0.086$ & $0.086$ & $0.859$ \\
    &  &  &  & 5000  & $\phantom{-}0.075$ & $0.041$ & $0.041$ & $0.561$ & $\phantom{-}0.077$ & $0.041$ & $0.041$ & $0.535$ \\    \hline
  \end{tabular}
\end{table}

\begin{table}[ht]
  \centering
  \caption{Simulation 3: Even if the instruments have the same distribution and all
    other assumptions are met, the TSIV estimators can still be biased
    if \Cref{assump:invariance-noise} (sampling homogeneity of noise) is
    violated. The reported numbers
    are obtained by 10000 realizations of data simulated from equations \eqref{eq:sim-iv},
    \eqref{eq:sim1-exposure}, \eqref{eq:sim1-outcome}, and
    \eqref{eq:sim1-uv}.}
  \label{tab:sim3}
  \begin{tabular}{|lllll|rrrr|rrrr|}
    \hline
    & & & & & \multicolumn{4}{c|}{TSTSLS} &
    \multicolumn{4}{c|}{Optimal TSIV} \\
    $\beta$ & $\rho^a$ & $\rho^b$ & $n^a$ & $n^b$ & Bias & SD & SE & Cover &
    Bias & SD & SE & Cover \\
    \hline
    1 & 0.5 & 0.5 & 1000 & 1000  & $-0.14$ & $0.279$ & $0.268$ & $0.899$ & $-0.14$ & $0.279$ & $0.268$ & $0.898$ \\
    &  &  &  & 5000  & $-0.13$ & $0.145$ & $0.131$ & $0.761$ & $-0.13$ & $0.145$ & $0.131$ & $0.762$ \\
    &  &  &  & 20000  & $-0.14$ & $0.097$ & $0.082$ & $0.555$ & $-0.14$ & $0.096$ & $0.082$ & $0.552$ \\
    &  &  & 5000 & 1000  & $-0.10$ & $0.294$ & $0.268$ & $0.892$ & $-0.10$ & $0.294$ & $0.268$ & $0.892$ \\
    &  &  &  & 5000  & $-0.09$ & $0.123$ & $0.122$ & $0.867$ & $-0.09$ & $0.123$ & $0.122$ & $0.867$ \\
    &  &  &  & 20000  & $-0.10$ & $0.067$ & $0.065$ & $0.667$ & $-0.10$ & $0.067$ & $0.065$ & $0.669$ \\
    &  &  & 20000 & 1000  & $-0.08$ & $0.271$ & $0.266$ & $0.939$ & $-0.08$ & $0.271$ & $0.266$ & $0.939$ \\
    &  &  &  & 5000  & $-0.08$ & $0.127$ & $0.120$ & $0.886$ & $-0.08$ & $0.127$ & $0.120$ & $0.886$ \\
    &  &  &  & 20000  & $-0.09$ & $0.062$ & $0.061$ & $0.707$ & $-0.09$ &
    $0.062$ & $0.061$ & $0.707$ \\
    \hline
  \end{tabular}
\end{table}

Even if $\mathbf{z}_i^a \overset{d}{=} \mathbf{z}_i^b$, the TSIV
estimators can still be biased if $v_i^a$ and $v_i^b$ have different
distributions and the instrument-exposure equation is not additive
(see the discussion after \Cref{assump:exposure-additive}). In our
third and final simulation, we generate the data from
equations \eqref{eq:sim-iv} and \eqref{eq:sim1-outcome} but replace
equations \eqref{eq:sim1-exposure} and \eqref{eq:sim1-uv} with
\begin{align}
  x_i^s &= I(0.2 \cdot (\mathbf{1}^T \mathbf{z}_i^s) + v_i^s > 0), \label{eq:sim3-exposure} \\
  (v_i^s,u_i^s) &\overset{i.i.d.}{\sim} \mathrm{N}\left(
    \begin{pmatrix}
      0 \\
      0 \\
    \end{pmatrix},
    \begin{pmatrix}
      \sigma_{vv}^s & \sigma_{uv}^s \\
      \sigma_{uv}^s & 1 \\
    \end{pmatrix}
  \right),~i=1,\dotsc,n^s,~s=a,b. \label{eq:sim3-uv}
\end{align}
In this simulation we use $\rho^a = \rho^b = 0.5$, $\sigma_{vv}^a =
1$, $\sigma_{vv}^b = 2$, and $\sigma_{uv}^s = 0.5
\sqrt{\sigma_{vv}^s}$, $s=a,b$.

The results of the third simulation are reported in
\Cref{tab:sim3}. Even though $\mathbf{z}_i^a \overset{d}{=}
\mathbf{z}_i^b$ in this simulation, the TSIV estimators are still
biased because the best linear approximations of the
instrument-exposure equation depend on the distributions of $v$, which
are different in the two samples.

\section{Application: The causal effect of body mass index on systolic
blood pressure}
\label{sec:real-data-analysis}

We apply the one-sample and two-sample IV methods to estimate the causal
effect of body mass index (BMI) on systolic blood pressure (SBP) using a real dataset
obtained from UK Biobank with 358,928 samples. As benchmarks, we first
apply ordinary least squares (OLS) and two IV methods (TSLS and LIML) to
the entire dataset with 407 correlated SNPs identified from a previous GWAS of BMI
\citep{locke2015genetic}. The results are reported in the first block
in \Cref{tab:ukbb}. The point estimate and confidence interval
obtained by OLS are much larger than those obtained by TSLS and LIML,
indicating there may be confounding in the observational
data. Unsurprisingly, the one-sample IV estimates agree with the
two-sample IV estimates using a random 50-50 split (second block in
\Cref{tab:ukbb}) and the
summary-data MR estimate reported in
\citet{zhao2018statistical} (third block in \Cref{tab:ukbb}).

Next we illustrate the performance of TSIV estimators with
heterogeneous samples. Because the UK Biobank population is mostly
homogeneous (most of samples are Europeans), we decide to subsample half
of the dataset in order to change the distribution of $9$ selected
SNPs. This artificially created subsample is then used as the exposure
(fourth block in \Cref{tab:ukbb}) or the outcome (fifth block in
\Cref{tab:ukbb}), while the other half of the dataset remains
unchanged and is used as the other sample in TSIV analyses. We find
that the TSIV point estimates using the two heterogeneous samples are
different from the benchmarks, though the differences are not
statistically significant due to increased standard error. Another
observation from \Cref{tab:ukbb} is that the TSTSLS estimator and the
optimal TSIV estimator always give very similar answers. This is not
surprising following the discussion in \Cref{rem:efficiency}.

\begin{table}
  \centering
  \caption{Results of the one-sample and two-sample IV analyses of the
    UK Biobank data.}
  \label{tab:ukbb}
  \begin{tabular}{|cc|c|c|c|c|}
    \hline
    Data & \# SNPs & Method & Estimate & Standard error & 95\% CI \\
    \hline
    \multirow{3}{*}{One-sample} & \multirow{3}{*}{407} & OLS & 0.7852 & 0.0068 & [0.7719,
    0.7985] \\
    & & TSLS & 0.4463 & 0.0366 & [0.3746, 0.5180] \\
    & & LIML & 0.3946 & 0.0392 & [0.3178, 0.4714] \\
    \hline
    \multirow{2}{*}{\begin{tabular}{l} Two-sample \\ (50-50 split) \end{tabular}} & \multirow{2}{*}{407} & TSTSLS & 0.4273 & 0.0514
    & [0.3266, 0.5280] \\
    & & Optimal\ TSIV & 0.4274 & 0.0514 & [0.3267, 0.5281] \\
    \hline
    \begin{tabular}{c} Two-sample \\ (summary data) \end{tabular} & 160
    & \begin{tabular}{c} MR-RAPS \\ \citep{zhao2018statistical} \end{tabular} &
    0.4017 & 0.1063 & [0.1934, 0.6100] \\
    \hline
    \multirow{2}{*}{\begin{tabular}{c} Two-sample \\ (subsampled
        exp.) \end{tabular}} &
    \multirow{2}{*}{9} & TSTSLS & 0.5199 & 0.1651 & [0.1963, 0.8435]
    \\
    & & Optimal\ TSIV & 0.5210 & 0.1651 & [0.1974, 0.8446] \\
    \hline
    \multirow{2}{*}{\begin{tabular}{c} Two-sample \\ (subsampled
        out.) \end{tabular}} &
    \multirow{2}{*}{9} & TSTSLS & 0.6500 & 0.1975 & [0.2629, 1.0371]
    \\
    & & Optimal\ TSIV & 0.6489 & 0.1975 & [0.2618, 1.0360] \\
    \hline
  \end{tabular}
\end{table}

\section{Summary and discussion}
\label{sec:summary-discussion}

In this paper we have derived a class of linear TSIV estimators when the two
samples are heterogeneous. Although the TSTSLS estimator is not
asymptotically efficient in general, it usually has great relative
efficiency and performs very similarly to the optimal TSIV estimator
in the numerical examples. Therefore there is little reason to abandon
the already widely-used TSTSLS in practice.

However, when trying to relax the linearity assumption, our
theoretical investigation suggests there are additional concerns about
using a two-sample IV analysis with heterogeneous samples.
\begin{enumerate}
\item Our (in fact any) TSIV analysis can only identify causal effect
  in the instrument-outcome sample (sample $b$). This is because we do
  not observe the outcome in sample $a$. This might limit the
  generalizability of the results of a real study.
\item Compared to the classical one-sample analysis, the TSIV analysis
  requires additional assumptions to link the two samples. One of the
  key assumptions is structural invariance (\Cref{assump:invariance}),
  which might be reasonable in some applications but unreasonable in
  others (especially if the two populations are drastically different).
\item Another important assumption in the two-sample setting is
  homogeneity of the distributions of the noise variables
  (\Cref{assump:invariance-noise}), which is necessary when the
  exposure equation is not additive. However, this assumption is
  untestable since we do not observe the exposure variable in one of
  the samples.
\item Unlike one-sample IV analysis, the heterogeneous two-sample IV
  analysis generally needs correct specification of the
  instrument-exposure equation.
\end{enumerate}
Our simulation examples show that violation of any of these three
requirements can lead to biased estimates and invalid statistical
inference. More real data examples are needed to evaluate the
importance of these concerns in practice.

The last point, that is the non-robustness of TSIV to model
misspecification and heterogeneous samples, is related to the notion
of ``invariant prediction''
\citep{peters2016causal}, ``autonomy''
\citep{haavelmo1944probability}, or ``stability''
\citep{pearl2009causality}. These notions are generally stronger as
they require invariance of the model under causal interventions. In
the problem considered in this paper, we require the exposure
predictions are invariant in the two heterogeneous samples. In this
view, the structural invariance (\Cref{assump:invariance}) is also not
necessary for the identification results. What's important is the
``predictive invariance'' in the two samples. In other words, even when
\Cref{assump:invariance} is violated so $f^a \ne f^b$, the causal
effect may still be identifiable if the best linear
approximations $\bm \gamma^a$ and $\bm \gamma^b$ defined in
\eqref{eq:gamma-projection} are the same.  We thank an associate editor
for pointing out this connection.


\bibliographystyle{chicago}
\bibliography{ref}

\end{document}